\documentclass[11pt]{amsart}
\PassOptionsToPackage{usenames,dvipsnames}{xcolor}
\usepackage{txfonts}
\PassOptionsToPackage{hyphens}{url}\usepackage{hyperref}
\usepackage{subfig,labelfig}
\usepackage{amscd,amssymb,mathabx}
\usepackage[arrow,matrix,graph,frame,poly,arc,tips]{xy}
\usepackage{graphicx,calrsfs}
\usepackage{mathtools}
\usepackage{tikz}
\usetikzlibrary{decorations.markings}
\usetikzlibrary{shapes}
\usetikzlibrary{arrows}
\usetikzlibrary{backgrounds}
\usetikzlibrary{calc}
\usepackage{mdwlist}
\usepackage{multirow}
\usepackage{enumerate}
\usepackage[shortlabels]{enumitem}

\usepackage{verbatim}
\usepackage{etoolbox}
\AtEndEnvironment{proof}{\setcounter{claim}{0}}

\DeclarePairedDelimiter{\form}{\langle}{\rangle}

\DeclarePairedDelimiter{\floor}{\lfloor}{\rfloor}

\newcommand\ba{\begin{align*}}
\newcommand\ea{\end{align*}}
\newcommand\be{\begin{enumerate}}
\newcommand\ee{\end{enumerate}}
\newcommand\bp{\begin{proof}}
\newcommand\ep{\end{proof}}
\newcommand\bpp{\begin{prop}}
\newcommand\epp{\end{prop}}
\newcommand\bpb{\begin{prob}}
\newcommand\epb{\end{prob}}
\newcommand\bd{\begin{defn}}
\newcommand\ed{\end{defn}}
\newcommand\bh{\begin{hint}}
\newcommand\eh{\end{hint}}

\newcommand\bN{\mathbb{N}}

\newcommand\bQ{\mathbb{Q}}

\newcommand\bZ{\mathbb{Z}}

\newcommand\bH{\mathbb{H}}

\newcommand\PP{\mathcal{P}}
\newcommand\TT{\mathcal{T}}

\DeclareMathOperator\Fibonacci{Fibonacci}
\DeclareMathOperator\im{im}
\DeclareMathOperator\Farey{Farey}
\DeclareMathOperator\hull{hull}

\DeclareMathOperator\Area{Area}

\newcommand\SL{\operatorname{SL}}

\newcommand\PSL{\operatorname{PSL}}

\DeclareMathOperator\cusp{cusp}

\newcommand\sse{\subseteq}
\newcommand\co{\colon}
\DeclareMathOperator\tr{tr}

\def\thetitle{{Optimal independent generating system for the congruence subgroups $\Gamma_0(p)$ and $\Gamma_0(p^2)$}}
\def\theauthors{{Nhat Minh Doan, Sang-hyun Kim, Mong Lung Lang, Ser Peow Tan}}
\usepackage{hyperref}
\hypersetup{
   colorlinks=false,
   plainpages,
   urlcolor=black,
   linkcolor=black
   pdftitle=   \thetitle,
   pdfauthor=  {\theauthors}
}

\theoremstyle{plain}
\newtheorem{thm}{Theorem}[section]
\newtheorem{lem}[thm]{Lemma}
\newtheorem{cor}[thm]{Corollary}
\newtheorem{prop}[thm]{Proposition}

\newtheorem*{principle*}{Principle}
\newtheorem*{claim*}{Claim}

\newtheorem{claim}{Claim}

\theoremstyle{remark}
\newtheorem{exmp}[thm]{Example}
\newtheorem{rem}[thm]{Remark}

\theoremstyle{definition}
\newtheorem{defn}[thm]{Definition}
\newtheorem{prob}{Problem}[section]

\begin{document}
\title[congruence subgroup]\thetitle
\date{\today}

\dedicatory{Dedicated to Professor Ravindra Kulkarni, on the occasion of his eightieth birthday.}

\author[Doan N. M.]{Nhat Minh Doan}
\address{Institute of Mathematics, Vietnam Academy of Science and Technology, Vietnam}
\email{dnminh@math.ac.vn}

\author[Kim S.-h.]{Sang-hyun Kim}
\address{School of Mathematics, Korea Institute for Advanced Study (KIAS), Seoul, Korea}
\email{skim.math@gmail.com}

\author[Lang M. L.]{Mong Lung Lang}
\address{Singapore}
\email{lang2to46@gmail.com}

\author[Tan S. P.]{Ser Peow Tan}
\address{Department of Mathematics, National University of Singapore, Singapore}
\email{mattansp@nus.edu.sg}

\begin{abstract}
Let $n$ be a prime or its square.
We prove that the congruence subgroup $\Gamma_0(n)$ admits a free product decomposition into cyclic factors in such a way that the $(2,1)$--component of each cyclic generator is either $n$ or $0$, answering a conjecture of Kulkarni.
We can also require that the Frobenius norm of each generator is less than $2n-1$.
A crucial observation is that if $P$ denotes the convex hull of
the extended Farey sequence of order $\floor{\sqrt{n}}$ in the hyperbolic plane $\bH^2$, then 
the projection $\pi\co \bH^2\to \bH^2/\Gamma_0(n)$ is injective on the interior of $P$ 
and 
each connected component of  $\pi(\bH^2)\setminus\pi(P)$ is either an order--three cone of area $\pi/3$ or an ideal triangle.
Denoting by $m(\Gamma_0(n))$  the minimum  of the largest denominator in the cusp set of $Q$ where $Q$ ranges over all possible special (fundamental) polygons for $\Gamma_0(n)$,
we establish the inequality
$ \floor{\sqrt{n}}\le m(\Gamma_0(n))\le \floor{\sqrt{4n/3}}$,
and completely characterize the cases in  which the bounds are achieved.
We also prove  analogous results when $n$ is the multiplication of two sufficiently close odd primes. 
\end{abstract}

\maketitle

\section{Introduction}\label{s:intro}
Let us denote the classical modular group as \[
\Gamma:=\PSL(2,\bZ)
=
\form{S,R\mid S^2=R^3=1},\]
where
\[
S=\begin{pmatrix}0&1\\ -1&0\end{pmatrix}\quad
\text{ and }\quad
R=\begin{pmatrix}0&1\\ -1&1\end{pmatrix}.\]
We also let
\[
T:=R^{-1}S=\begin{pmatrix}1&1\\ 0&1\end{pmatrix}.\]
Each subgroup $\Lambda$ of $\Gamma$ admits a free product decomposition
\[
\Lambda=\bigast_i \form{s_i},\]
where the order of each $s_i$ is either two, three, or infinite. Such a generating set $\{s_i\}$ is called an \emph{independent generating system} (or \emph{independent system}) for $\Lambda$. When $[\Gamma:\Lambda]<\infty$, the Reidemeister--Schreier process yields a presentation of $\Lambda$, and a suitable simplification of the presentation would yield an independent system for $\Lambda$. 
In the case of the congruence subgroup
\[
\Gamma_0(n):=
\left\{A\in\Gamma\co A\equiv \begin{pmatrix} * & * \\ 0 & *\end{pmatrix} \pmod{n}
\right\},\]
such a process was given by Rademacher~\cite{Rademacher1929} for a prime $n$, and by Chuman~\cite{Chuman1973} for a general integer $n$.

Kulkarni~\cite{Kulkarni1991} gave a geometric and much more practical algorithm to find an independent system for every finite index subgroup $\Lambda$ of $\Gamma$.
In general, a fundamental domain of $\Lambda$ equipped with a side-pairing rule yields a generating set for $\Lambda$, but this is not guaranteed to be independent. 
Kulkarni utilized/constructed  a certain particularly useful type of fundamental domain (\emph{special polygon}) $P^*$. He encoded this information as a \emph{Farey symbol}, which is equivalent to a maximal ideal polygon $P$ inside $P^*$ with each side of $P$ suitably labeled.
The independent system in the \emph{matrix form} can be read off immediately from the output of his algorithm. See Section~\ref{s:special}.

For $\Gamma_0(n)$, Kulkarni's algorithm relies on congruence relations among the denominators of certain rational numbers, which correspond to the cusps of a special polygon. After long and meticulous computations, he noted~\cite{Kulkarni1991} that for $\Gamma_0(p)$ with $p$ prime: 
\begin{quote}
``For $p < 100$ we observed that the congruences which need to be satisfied for constructing these Farey symbols can actually be lifted to equalities in natural numbers''
\end{quote}
and suggested that this could be true in general. Here, the ``congruences'' mean that the $(2,1)$ components of the resulting independent generators are $0$ modulo $p$. Our main theorem below  establishes his conjecture.
\begin{thm}\label{t:main}
If $n$ is a prime or the square of a prime, then $\Gamma_0(n)$ admits 
a free product decomposition
\[
\Gamma_0(n) =\bigast_{i=1}^k\form*{g_i}\bigast \form{T }\]
for some  matrices $g_i$ such that $(g_i)_{2,1}=n$;
furthermore, we can require that $|\tr g_i|\le n-2$ 
and $\|g_i\|_F<2n-1$.
\end{thm}
In the above $\|g\|_F:=\sqrt{\tr(gg^t)}$ denotes the Frobenius norm of a matrix $g$.
We also establish an analogous result for the case when $n=pq$ for some sufficiently close odd primes $p$ and $q$ (Theorem~\ref{thm:twin}). 
For a general integer $n$,
it was recently announced~\cite{Steiner2022} that 
$\Gamma_0(n)$ admits a (not necessarily independent) generating set consisting of matrices with Frobenius norms at most $O_\epsilon(n^{1+\epsilon})$; see also~\cite{MR1146784} for the prime power case.
If $n$ is a prime or its square then the number of the free factors $k$ is bounded above by a fixed linear function on $n$, as can be seen from formulae in Section~\ref{ss:gamma0}.




Recall the \emph{Farey sequence $\Farey_k$ of order $k$} is the set of rational numbers in $[0,1]$ whose denominators are at most $k$. 
Let us write \[\Farey_k^*:=\Farey_k\cup\{\infty\}\sse\bQ\cup\{\infty\}=:\hat\bQ\sse S^1=\partial\bH^2,\]
which is called the \emph{extended Farey sequence of order $k$}.
A special polygon for $\Lambda\le \Gamma$ is said to be \emph{optimal} if the maximum value of the denominators of its cusps realizes the minimum value among all possible special polygons for $\Lambda$, and we denote by $m(\Lambda)$ this value. Note that $\Lambda$ may have several different optimal special polygons. We prove Theorem~\ref{t:main} in the process of establishing an  efficient bound on the value $m(\Gamma_0(n))$.

\begin{thm}\label{t:bound}
The following holds.
\be
\item
For all integers $n\ge2$, we have that
$m(\Gamma_0(n))\ge \floor{\sqrt{n}}$;
moreover,  equality holds iff $n$ belongs to the list
\[n\in\{2,3,4, 5,7,9, 11, 13, 17, 19,25, 29,31, 37, 49, 53, 67, 83, 127,173\}.\]
\item
Assume $n=p$ or $p^2$ for a prime $p$.
Then 
$m(\Gamma_0(n))\le \floor*{\sqrt{4n/3}}$,
and the equality holds iff $n$ is of the form $sa+tb$
for some $ s,t,a,b \in \bN_+$ satisfying 
\begin{itemize}
   \item $s+t>a>t\geq b \geq a-s$.
    \item $3(a+1)^2>4(sa+tb)$.
\end{itemize}
\ee
\end{thm}

\begin{rem}
We note that the list in part (1) of Theorem~\ref{t:bound} happens to consist of primes and prime squares only.
\end{rem}

Let $v:=\floor{\sqrt{n}}$, and
 let $P_v$ be the convex hull of $\Farey_v^*$ in the hyperbolic plane.
A starting point towards the above theorem is the observation that 
for an arbitrary integer $n\ge 2$, 
the projection \[\pi\co\bH^2\to \bH^2/\Gamma_0(n)\] maps the interior of $P_v$ injectively;
see Lemma~\ref{lem:extends}.
This readily establishes an effective lower bound of $v$ in the above theorem.
The equality condition will follow from an asymptotic estimate of the number of ideal triangles in $P_v$.

The second key result we show is that 
under the extra hypothesis that $n=p$ or $p^2$ for a prime $p$,
each connected component of 
\[
\pi(\bH^2)\setminus\pi(P_v)\]
is either an order--three cone
or an ideal triangle; see
 Remark~\ref{rem:estimate}. Here, 
an \emph{order--three cone} means the image of a triangle with angles $(0,0,2\pi/3)$ with the two rays at the angle $2\pi/3$ glued.
We also have a geometric proof of the inequality
\[
\floor*{\frac{n(1+1/p)}3}\ge\sum_{1\le i\le \floor{\sqrt{n}}}\phi(i),\]
by showing that the difference is precisely the number of such ideal triangles (constituting connected components).
The equality holds iff $n$ belongs to the list in Theorem~\ref{t:bound} (1).
As a consequence, this difference can be expressed in terms of the number of solutions to certain diophantine equations; see Remark~\ref{rem:estimate}.

In the process of proving this result, we show that every side-pairing not fixing $\infty$ can be realized by a matrix whose $(2,1)$-entry is $n$, or equivalently, all the congruences satisfied by the denominators can be lifted to equalities, thereby establishing Theorem \ref{t:main}. 

Proving the optimality of our upper bound $\floor*{\sqrt{4n/3}}$ 
at least in the case of a prime or its square
requires more work and is related to the expressibility of integers as certain quadratic forms in one or two variables. This is done in Section \ref{s:optimal}. 
We remark that with computer verification one can make a similar assertion without the hypothesis $n\ge37$.
\begin{prop}
Let $n\ge37$ be a prime or its square.
Then we have \[m(\Gamma_0(n))=\floor*{\sqrt{4n/3}}\] if and only if there exist $s,t,a,b\in\bN_+$ such that
\begin{equation}
\begin{cases}
&n = sa+tb,\\
    &s+t>a=\floor{\sqrt{4n/3}}>t\geq b \geq a-s\ge1.
\end{cases}
\end{equation}
\end{prop}

The following corollary is a special case by choosing $a=s+t-2$ and $b = t-2$,
and by further choosing  that $t=s+1$.
\begin{cor}\label{cor:upperbound}
Let $p\ge 37$ be a prime in  one of the following forms:
\begin{itemize}
\item $p=s^2+t^2+st-2s-2t,$\\
for some $s,t\in \bN_+$ with $3 \leq s<t<s+1+2\sqrt{s+1}$;
\item $p=3s^2-s-1,$\\
for some $s\in \bN_+$.
\end{itemize}
Then we have that $ m(\Gamma_0(p))=\floor*{\sqrt{4p/3}}$.
\end{cor}
Iwaniec \cite{Iwaniec1974} proved that there are infinitely many primes of the first form above, but we do not know if they satisfy the extra conditions on the bounds of $s$ and $t$.
The second condition above in particular implies that the upper bound $\floor*{\sqrt{4p/3}}$ can be achieved for infinitely many primes $p$ if we \emph{assume} a particular case of the Bunyakovsky conjecture~\cite{Bunyakovsky1857}: there exist infinitely many primes $p$ of the form $3s^2-s-1$ with $s\in\bN_+$.
\medskip

The very tight connection between the geometry and the arithmetic for the groups $\Gamma_0(n)$ established by the above results was somewhat surprising and hints at a  deeper and more mysterious relation for congruence subgroups. 
For a general subgroup $\Lambda\le \Gamma$ of index $n$, one cannot hope to obtain this type of bound,  see Example \ref{ex:Fibonacci} where $m(\Lambda)$ is exponential in $n$. The upper bound established above does depend on $n$ being a prime or its square. For example, in $\Gamma_0(8)$ the upper bound above does not hold; see Example \ref{ex:8}.

\section{Preliminaries}\label{s:special}
In this section, we give a brief exposition on Kulkarni's construction of \emph{special polygons} for congruence subgroups in $\Gamma$.


\subsection{Special polygons}\label{ss:special}
The group $\Gamma:=\PSL(2,\bZ)$ acts faithfully on the upper half-plane $\bH^2$ by orientation-preserving isometries,
so that  each element of $\Gamma$ is identified with the corresponding linear fractional transformation. 
We let $S$ be the rotation of angle $\pi$ about $\sqrt{-1}$,
and $R$ be the rotation of angle $2\pi/3$ about $e^{\pi i/3}$.
We have a group presentation
\[\Gamma=\form{S,R\mid S^2=R^3=1}.\]
We set $T(z):=R^{-1}S(z)=z+1$.
This group $\Gamma$ leaves invariant the circularly ordered set \[\hat\bQ=\bQ\cup\{\infty\}\sse\partial\bH^2\approx S^1.\]
Unless specified otherwise, we always assume that an expression $a/b$ for an element of $\hat\bQ$ is in the \emph{reduced form}, that is $\gcd(a,b)=1$ and $b\ge0$; such an expression is unique except for
\[
\infty = 1/0 = (-1)/0 =-\infty.\]
For distinct $x,y\in\hat\bQ$ the ordered pair $(x,y)$ is often identified with the oriented infinite geodesic line in $\bH^2$ moving from $x$ to $y$.

A positively oriented triangle in the $\Gamma$--orbit of the geodesic triangle $\Delta_0:=(0,1,\infty)$ is called a \emph{Farey triangle}. 
The \emph{Farey tessellation} $\TT$ is the 2--complex structure on $\bH^2\cup\hat\bQ$ having the Farey triangles as the 2--cells.
Each triangle in the $\Gamma$--orbit of 
 $\Delta_0^*:=(0,e^{\pi i/3},\infty)$ is called a \emph{special triangle}.
We let $\TT^*$ denote the subdivision of $\TT$ consisting of special triangles so that each Farey triangle is subdivided into three special triangles.
Each special triangle is a fundamental domain for $\Gamma$.

The pair of the endpoints of a geodesic line in the Farey tessellation is called  a \emph{Farey pair}.
For a Farey pair $(x,y)$ we define their \emph{mediant} $x\oplus y$ as the unique vertex in $\hat \bQ$ such that $(x,x\oplus y, y)$ is a (positively oriented, as usual) Farey triangle.
Algebraically, $(a'/a,b'/b)\in\hat\bQ\times\hat\bQ$ 
is a Farey pair if $ab'-a'b=\pm1$;
if we further assume $a'/a,b'/b\in\bQ$, then we see that
\[
\frac{a'}a\oplus \frac{b'}b = \frac{a'+b'}{a+b}.\]

For a geodesic polygon $P\sse\bH^2$ we let $E(P)$ and $\cusp(P)$ denote the sets of the positively oriented sides and of the ideal vertices of $P$, respectively.
For a side $e$, we let $\bar e$ denote its reversely oriented geodesic.

\bd\label{d:polygon}
Let $\Lambda$ be a subgroup of the classical modular group $\Gamma$.
\be
\item A \emph{normalized ideal polygon} is a finite convex subcomplex $P$ of $\TT$
such that $\infty$ and $0$ are vertices of $P$;
in the case when $\Lambda=\Gamma_0(n)$ for some $n$, we further require that $(\infty,0)$ is a positively oriented side of $P$;
\item 
A normalized ideal polygon
whose interior maps injectively under the quotient $\bH^2\to\bH^2/\Lambda$
is called a \emph{$\Lambda$--injective ideal polygon}, or simply a \emph{$\Lambda$--polygon}.
Such a polygon is \emph{maximal} if it is maximal under inclusion.
\item\label{p:special} A \emph{special polygon $P^*$ for $\Lambda$} is a convex subcomplex of $\TT^*$ which is a fundamental domain of $\Lambda$
such that the convex hull of $\cusp P^*$ is normalized.
In this case, the set
\[\{g\in\Lambda\mid P^*\cap gP^*\text{ is a side of }P^*\}\]
is called the \emph{side-pairing rule} corresponding to $P^*$.
\ee
\ed

We observe that the convex hull of $\cusp(P^*)$ in part (\ref{p:special}) is a maximal $\Lambda$--polygon.
Note also that a side-pairing rule is closed under taking inverses.
A remarkable discovery of Kulkarni in~\cite{Kulkarni1991} is that a special polygon provides a geometric tool for the purely algebraic problem of efficiently finding an independent system for $\Lambda\le \Gamma$.

\begin{thm}[\cite{Kulkarni1991}]\label{thm:kulkarni}
Every finite index subgroup $\Lambda$ of $\Gamma$ admits a special polygon, in such a way that if $U$ is a maximal subset of the side pairing rule not containing two distinct elements which are inverses to each other
then $U$ is an independent system (in Rademacher's sense) for $\Lambda$.
\end{thm}

A key idea in Kulkarni's construction is that every $\Lambda$--polygon $P$ can be enlarged to some maximal $\Lambda$--polygon, and then to a special polygon for $\Lambda$; this follows from the fact that each side of such a $P$ can be classified as below.

\begin{lem}\label{lem:types}
Let $P$ be a $\Lambda$--polygon for some 
$\Lambda\le\Gamma$.
Suppose $e$ is a side of $P$, and denote by $\Delta$ the unique Farey triangle that shares only the side $e$ with $P$. 
Then  exactly one of the following holds.
\be[(i)]
\item $e$ is \emph{even}, meaning  a unique $g_e\in \Lambda$ reverses the orientation of $e$.
\item $e$ is \emph{odd}, meaning  a unique $g_e\in\Lambda$ rotates $\Delta$ by the angle $2\pi/3$;
\item $e$ is \emph{paired} to a different side $f$ of $P$, meaning  a unique $g_e\in\Lambda$ 
satisfies $\bar f= g_e(e)$;
\item $e$ is \emph{free}, meaning  $P\cup\Delta$ is a $\Lambda$--polygon. \ee
\end{lem}

In particular, we have a \emph{labeling map} \[\sigma\co E(P)\to \{-4,-3,-2\}\cup\bN_+\] so that
the preimages of $-4,-3,-2$ are free, odd and even respectively, 
and the preimage of each positive integer is either empty or consists of two paired sides. 
We slightly abuse the terminology and say the set \[\{g_e\in\Lambda\mid e\text{ is a non-free side of }P\}\] is the \emph{side-pairing rule} for $P$.
If $P$ has no free sides, then $P$ is maximal and included in some special polygon $P^*$ with the same cusp set.
Hence, the above set coincides with the side-pairing rule of $P^*$ in the sense that was defined earlier. 
In the case when $P$ is a maximal $\Lambda$--polygon with a labeling map $\sigma$, then the pair
\[
(\cusp P,\sigma)\]
is often called a \emph{Farey symbol} in the literature; see~\cite{Kulkarni1991}.

To construct a special polygon for a given group $\Lambda$, one  starts with the Farey triangle $(0,1,\infty)$ and continues adding a Farey triangle to a remaining free side to obtain a $\Lambda$--polygon.
If $\Lambda$ is of finite index,
this process must terminate at a maximal $\Lambda$--polygon $P$.
One then finally adds a special triangle to each odd side of  $P$, and obtains a special polygon as required in Theorem~\ref{thm:kulkarni}.

Conversely to Theorem~\ref{thm:kulkarni}, Kulkarni also noticed that if $P$ is a finite area convex subcomplex of $\TT$ with each side labeled by even, odd or a positive integer in such a way that each positive integer is assigned to either zero or two edges (called paired edges) then there uniquely exists a finite index subgroup $\Lambda_P\le\Gamma$ such that the edges of $P$ are classified as in the above lemma with respect to $\Lambda_P$, and such that $\Lambda_P$ admits a special polygon $P^*$ with $\cusp(P^*)=P$.
In particular, $\Lambda_P$ is generated by the side-pairing rule.

\begin{exmp}\label{ex:Fibonacci}
Let us fix a positive integer $n$. We define a  sequence \[
r_0:=0<r_2<r_4<\cdots<r_5<r_3<r_1:=1\]
by the condition that 
$r_{i+2}=r_{i+1}\oplus r_i$
for $i=0,\ldots, n-2$. In other words, we have $r_i:=\Fibonacci_{i}/\Fibonacci_{i+1}$.
We let $P(n)$ be the ideal polygon with the cusp set $\{r_i\}_{0\le i\le n}\cup\{\infty\}$,
and label the sides by the following rule.
\[
\sigma(e):=
\begin{cases}
1,&\text{ if }e=(-\infty,0)\text{ or }e=(1,\infty),\\
2,&\text{ if }e=(r_n,r_{n-1})\text{ or }e=(r_n,r_{n-2})\text{ up to orientation},\\
-2,&\text{ otherwise}.
\end{cases}\]
Then there exists a finite index subgroup $\Lambda_{P(n)} \le \Gamma$ generated by the side-pairing rule of $P(n)$.
The group $\Lambda_{P(n)}$ admits a special polygon,  which coincides with $P(n)$ since there are no odd sides.
The vertices created converge to $1/\phi$ where $\phi$ is the golden ratio.
The corresponding index satisfies 
\[
[\Gamma\co \Lambda_{P(n)}]/3=n,\]
which is the number of the Farey triangles in $P(n)$ after the final step. As there was essentially only one choice at each step of the construction once the initial ideal triangle $(0,1,\infty)$ is chosen, this is the unique special polygon for $\Lambda_{P(n)}$ up to the action of $\Gamma$.
We record for a later reference that the largest denominator in $\cusp(P(n))$ is exponential in $n$.
\end{exmp}

\subsection{The congruence subgroup $\Gamma_0(n)$}\label{ss:gamma0}
We are mostly concerned with the case  $\Lambda=\Gamma_0(n)$.
In the hyperbolic 2--orbifold
\[\bH^2/\Gamma_0(n),\]
we denote by $v_\infty(n),v_2(n)$ and $v_3(n)$ the numbers of cusps, order--two cone points and order--three cone points respectively.
It is well known that
\[
[\Gamma:\Gamma_0(n)]=n\prod_{p\divides n} \left(1+1/p\right).
\]
Let $\PP(n)$ denote the set of prime divisors of $n$.
We have~\cite[Page 25]{Shimura1971} that
\[
v_\infty(n)=\sum_{d\divides n} \phi\left(\gcd(d,n/d)\right),\]
and that 
\begin{align*}
v_2(n)&=\begin{cases} 0, &\text{ if }4\divides n\text{ or if }\PP(n)\cap(4\bZ-1)\ne\varnothing,\\
2^k, &\text{otherwise, with }k:=\#\left(\PP(n)\cap (4\bZ+1)\right).
\end{cases}\\
v_3(n)&=\begin{cases} 0, &\text{ if }9\divides n\text{ or if }\PP(n)\cap(3\bZ-1)\ne\varnothing,\\
2^k, &\text{otherwise, with }k:=\#\left(\PP(n)\cap (3\bZ+1)\right).
\end{cases}
\end{align*}
The genus $g$ can be computed from
the Riemann--Hurwicz formula. Namely,
\[
\frac\pi3[\Gamma:\Gamma_0(n)]=\Area(\bH^2/\Gamma_0(n))
=-2\pi\left(2 - 2g -v_\infty(n) - \frac{2v_3(n)}3 - \frac{v_2(n)}2\right).\]
The number of each free factor in the following decomposition 
is computed:
\[
\Gamma_0(n)\cong \left(\bigast_{i=1}^{v_2(n)} \bZ/2\bZ\right)\bigast
\left(\bigast_{i=1}^{v_3(n)} \bZ/3\bZ\right)\bigast
\left(\bigast_{i=1}^{2g+v_\infty(n)-1} \bZ\right).\]

One particular advantage of the approach in Theorem~\ref{thm:kulkarni} is that matrix representatives of an independent system can be read off directly during the process since the choice of the element $g_e$ is straightforward from the label and the endpoints. Let us record one consequence that will be relevant to our purpose. Recall our convention that rational numbers are assumed to be written in reduced form.
\begin{lem}[\cite{Kulkarni1991,CLLT1993}]\label{lem:2-1-compo}
Let $P$ be a $\Gamma_0(n)$--polygon,
and let $e=(a'/a,b'/b)$ be a side of $P$ not incident at $\infty$. 
Then the following hold.
\be
\item The side $e$ is even iff $n\divides a^2+b^2$;
\item The side $e$ is odd iff $n\divides a^2+ab+b^2$;
\item The side $e$ is paired to another side $f=(c'/c,d'/d)$ iff $n\divides ac+bd$.
\ee
\end{lem}


Suppose $P$ is a $\Gamma_0(n)$--polygon, whose cusp set is
written as
\[
x_{-1}=-\infty<x_0=0<x_1<\cdots<x_\ell=1<\infty=x_{-1}.\]
In particular, the two sides $(-\infty,0)$ and $(1,\infty)$ are paired.
The number $\ell$ of ideal triangles in a special polygon $P^*$ for $\Gamma_0(n)$ is given by 
\[ u(n):=\frac{\Area(P^*)-\pi v_3(n)/3}\pi =\frac{[\Gamma\co\Gamma_0(n)]- v_3(n)}3=\#\cusp (P^*)-2.\]
In particular, if $n$ is a nontrivial power of a prime $p$ then $u(n)=\floor*{n(1+1/p)/3}$.
The labeling map $\sigma$ of $P$ is often succinctly encoded as a tuple
\[(\sigma(-\infty,x_0=0):=1, \sigma(x_0,x_1),\ldots,\sigma(x_\ell=1,\infty):=1).\]

In general, the \emph{denominator sequence} of a sequence $A=\{b_i'/b_i\}$ in $\hat\bQ$ is defined as
\[
d(A):=d  A=\{b_i\}.\]
Let us now assume the denominator sequence of $\cusp P$ is written as
\[d (\cusp P):=\{a_{-1}=0,a_0=1,a_1,\ldots,a_\ell=1\}.\]
We call each consecutive pair $(a_i,a_{i+1})$  a \emph{side} of $d (\cusp P)$, with indices taken cyclically.  
Such a side $(a_i,a_{i+1})$ 
inherits its labeling from the side $(x_i,x_{i+1})\in E(P)$.

A $\Gamma_0(n)$--polygon $P$ is uniquely determined by the denominator sequence of $\cusp(P)$; this is due to the following general observation.
\begin{lem}\label{lem:side}
If $(a,b)$ is a pair of positive co-prime integers,
then there uniquely exist nonnegative integers $a',b'$ such that $a'/a$ and $b'/b$ are in reduced forms comprising a Farey pair, and such that
\[
0\le a'/a<b'/b\le 1.\]
\end{lem}
\bp
The equation $ax-by=1$ admits a unique  solution $(x,y)$ satisfying
$0\le y<a$ and $1\le x\le b$. This is equivalent to the condition
that $y/a$ and $x/b$ are a Farey pair in $[0,1]$, as required.
\ep

\begin{figure}[htb!]
\leavevmode \SetLabels
\L(.18*.15) $\frac{0}{1}$\\%
\L(.33*.15) $\frac{1}{4}$\\%
\L(.39*.15) $\frac{1}{3}$\\%
\L(.49*.15) $\frac{1}{2}$\\%
\L(.8*.15) $\frac{1}{1}$\\%
\endSetLabels
\begin{center}
\AffixLabels{\centerline{\includegraphics[width=10cm,angle=0,scale=.8,trim={0 0 0 .6cm},clip]{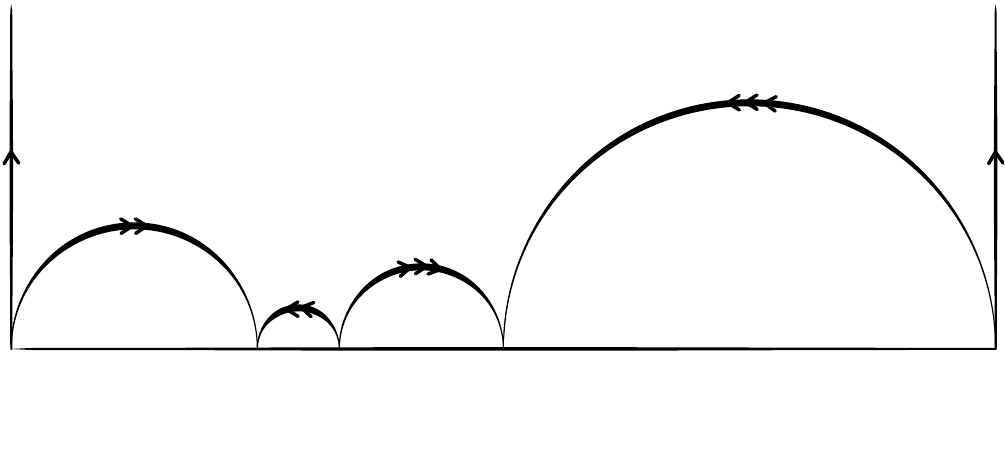}}}
\vspace{-24pt}
\end{center}
\caption{One of the two possible special polygons for $\Gamma_0(8)$. The denominator sequence is $\{ 0, 1,4, 3,2,1\}$.} 
\label{sp8}
\end{figure}

From the two preceding lemmas, we see that 
 it suffices to work only with the {denominator sequence} 
when finding and enlarging a $\Gamma_0(n)$--polygon $P$.
\begin{exmp}\label{ex:8}  The denominator sequence
\[ d (\cusp P)=\{0, 1, 4,3,2, 1\}\]
uniquely determines a normalized ideal polygon $P\le\TT$ as
\[ \cusp(P)=  \{-\infty ,0, 1/4, 1/3,1/2,1\}.\] 
In the case when $n=8$, 
it is easy to deduce from Lemma~\ref{lem:2-1-compo} and from the denominator sequence that
the labeling map of $P$ is given as
\[
\sigma
=(1,2,2,3,3,1).
\]
See Figure~\ref{sp8}. 
Since there are no odd or free sides, we see $P$ is a special polygon for $\Gamma_0(8)$.
There exists a different special polygon for $\Gamma_0(8)$. Namely, one can verify that $d (\cusp P')=\{0,1,2,3,4,1\}$ and $\sigma'=(1,3,3,2,2,1)$ also determine a special polygon $P'$ for the same group.
\end{exmp}

\begin{exmp}\label{ex:smallp}
Let us work out all possible denominator sequences and the labelings of maximal ideal polygons $P$ for $\Gamma_0(n)$ with $n = 2,3,5$ and $7$.
\begin{itemize}
    \item $n=2$: $d (\cusp P)=\{0,1,1\}$, $\sigma=(1,-2,1)$.
    \item $n=3$: $d (\cusp P)=\{0,1,1\}$, $\sigma=(1,-3,1)$.
    \item  $n=5$: $d (\cusp P)=\{0,1,2,1\}$, $\sigma=(1,-2,-2,1)$.
    \item $n=7$: $d (\cusp P)=\{0,1,2,1\}$, $\sigma=(1,-3,-3,1)$.
\end{itemize}
\end{exmp}
Note that $n=2,3$ and $n=5,7$ share the same cusps with different labelings.

\section{Bounds on the denominators of cusps}\label{s:order}
As briefly discussed in Section~\ref{s:intro},
a key step in our proof for Theorem~\ref{t:main} is to establish effective lower and upper bounds for the following quantity:
\begin{align*}
m(\Gamma_0(n))
:= &\min\left\{
\max  d (\cusp P) \mid P\text{ is a maximal }\Gamma_0(n)\text{--polygon}\right\}\\
= &\min\left\{
\max  d (\cusp P^*) \mid P^*\text{ is a special polygon for }\Gamma_0(n)\right\}
\end{align*}
Recall that a special polygon $P^*$ for $\Gamma_0(n)$ is \emph{optimal}
if the maximum of $d (\cusp(P^*))$ realizes such a minimum.
We will actually deduce Theorems~\ref{t:main} and~\ref{t:bound}
from stronger results, 
which are Theorem~\ref{t:bound}
and the following.

\begin{thm}\label{thm:bound}
If $n$ is a prime or the square of a prime,
and if $v:=\floor{\sqrt{n}}$,
then a maximal $\Gamma_0(n)$--polygon $P$ can be chosen so that
  the following holds with $U$ being a maximal subset of the side-pairing rule of $P$ not containing two distinct elements that are inverses to each other (cf.~Theorem~\ref{thm:kulkarni}).
\begin{itemize}
\item $\cusp P =\Farey_v^*\sqcup \{c_1<\cdots<c_k\}$ such that  for all $i<j$ the points $c_i$ and $c_j$ are not Farey neighbors;
\item Each $c_i$ and the two adjacent points of $c_i$ in $\cusp P$ form a Farey triangle $T_i$, and each of the two sides incident on $c_i$ in $T_i$ is paired to some side of the convex hull of $\Farey_v^*$;
\item For each element $g\in  U$ not fixing $\infty$, the $(2,1)$ component of $g$ is exactly $n$. Furthermore, 
$|\tr g|\le n-2$ and  $\|g\|_F:=\sqrt{\tr (gg^t)}<2n-1$.
\end{itemize}
\end{thm}

\begin{figure}[htb!]
\leavevmode \SetLabels
\L(.1*.15) $\frac{0}{1}$\\%
\L(.295*.15) $\frac{1}{4}$\\%
\L(.18*.58) \small even\\%
\L(.75*.58) \small even\\%
\L(.36*.15) $\frac{1}{3}$\\%
\L(.487*.15) $\frac{1}{2}$\\%
\L(.48*.6) \large $P$\\%
\L(.618*.15) $\frac{2}{3}$\\%
\L(.68*.15) $\frac{3}{4}$\\%
\L(.875*.15) $\frac{1}{1}$\\%
\endSetLabels
\begin{center}
\AffixLabels{\centerline{\includegraphics[width=10cm,angle=0,trim={0 0 0 1.1cm},clip]{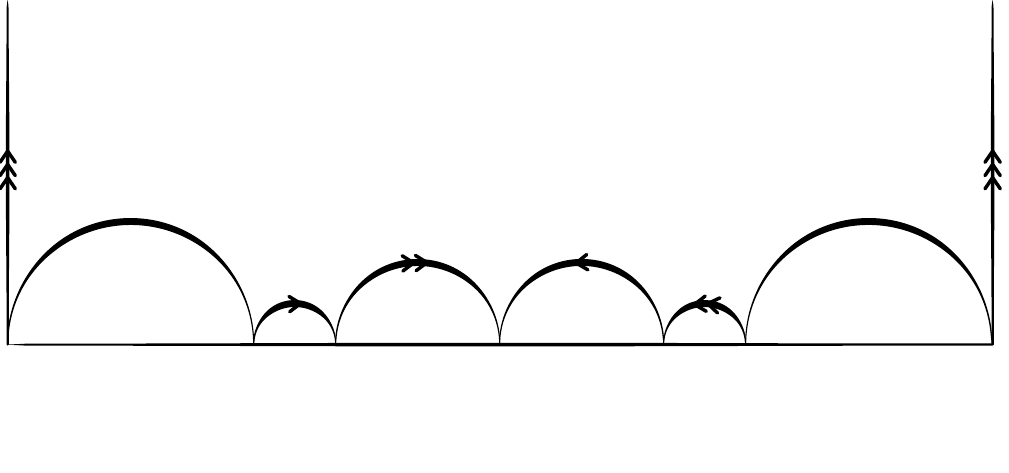}}}
\vspace{-24pt}
\end{center}
\caption{$m(\Gamma_0(17))=4$.} 
\label{fig:sp17}
\end{figure}

\begin{exmp}
When $p=17$, we have that
\[
\floor*{\sqrt{p}}=\floor*{\sqrt{4p/3}}=4.\]
In Figure~\ref{fig:sp17}, we illustrate a maximal $\Gamma_0(17)$--polygon $P$ and its side-pairing rule such that $\cusp(P)=\Farey_4^*$. This polygon $P$ happens to be a special polygon for $\Gamma_0(17)$ since there are no free or odd sides.
\end{exmp}

In the first subsection, we prove the former part of Theorem~\ref{t:bound}.
Theorem~\ref{thm:bound} and the latter part of Theorem~\ref{t:bound} is considered in the second subsection. 
\emph{Throughout this section, we let $n\ge2$ be an integer, and set
\[v:=\floor{\sqrt{n}},\qquad dF:=d \Farey_v^*.\]
We also let $P_v$ denote the convex hull of $\Farey_v^*$.}

\subsection{Establishing the lower bound}\label{ss:algorithm}
Recall we have denoted by $u(n)$ the number of ideal triangles in a maximal $\Gamma_0(n)$--polygon.
The \emph{totient summatory function} is defined as
\[
\Phi(n):=\sum_{i=1}^n\phi(i),\]
which coincides with the number of reduced fractions in $[0,1)$ with denominators up to $n$. 
Let us note the following preliminary estimate.
\begin{lem}\label{lem:lower1}
We have that
\[
\min\left\{w\mid  \Phi(w)\ge u(n)\right\}
\le
m(\Gamma_0(n))
\le
((1+\sqrt{5})/2)^{u(n)}.
\]
Furthermore, we have that
\[
\liminf_{n\to\infty} \frac{m(\Gamma_0(n))}{[\Gamma\co\Gamma_0(n)]^{1/2}}\ge \frac\pi3.\]
\end{lem}
\bp
Fix an integer $n\ge2$ and set $m:=m(\Gamma_0(n))$.
Since there exists a  maximal $\Gamma_0(n)$--polygon $Q$ whose cusp set is contained in $\Farey_m^*$,
we have
\[
u(n)= \#\cusp(Q)-2
\le \#\Farey_m^*-2=\Phi(m).\]
This gives the first inequality.
For the second one, it suffices to note that
each cusp of $Q$ belongs to a Farey triangle that can be joined to the Farey triangle $(0,1,\infty)$ by a sequence of adjacent Farey triangles,
and that such a sequence has length at most $u(n)$. See also Example~\ref{ex:Fibonacci}.

To see the last part of the lemma, let us note the classical fact
 \[\Phi(m)=
 {3}m^2/\pi^2 + O (m(\log m)^{2/3}
 (\log\log m)^{4/3}).\]
See~\cite{Walfisz1963}.
For $n\ge4$,  we claim that the (crude) bound
\[
v_3(n)\le \sqrt{n}\]
holds. 
Indeed, suppose $n=p_1^{e_1}\cdots p_{k}^{e_{k}}$ is the prime factorization with \[p_1<\cdots<p_{k}.\]
We only consider the case that $k\ge2$ and $p_i\in\{3\}\cup(3\bZ+1)$ for all $i$, since the claim is trivial otherwise.
We observe
\[
p_1 \cdots p_k
\ge 3\cdot\prod_{i=1}^{k-1}(6i+1)\ge 4^k.\]
It follows that 
$v_3(n)=2^k\le\sqrt{p_1 \cdots p_k}\le\sqrt{n}$.

Assuming $n$ is sufficiently large, we have that
\[
[\Gamma\co\Gamma_0(n)]-\sqrt{n}
\le 3u(n)\le 3\Phi(m)\le
9m^2/\pi^2+3m\log m.\]
Since $[\Gamma\co\Gamma_0(n)]>n$, we obtain the latter part of the conclusion.
\ep

The following elementary observation will be useful for us.
\begin{lem}\label{lem:elem}
If $a,b,x,y$ are integers in $(0,v]$ such that $x+y\le v$
and such that $\gcd(a,b)=1$, then we have that
\[
n\not\divides ax+by,\quad n\not\divides x^2+y^2,\quad n\not\divides x^2+xy+y^2.\]
\end{lem}
\bp
Since $a$ and $b$ are co-prime, we have  $\min(a,b)<v$.
It follows that
\[
2\le ax+by <v(x+y)\le n.\]
We also have that $x^2+y^2<x^2+xy+y^2 <(x+y)^2\le n$, implying the conclusion.\ep

A parabolic element in $\Gamma$ is a \emph{translation} iff it fixes $\infty$.
\begin{lem}\label{lem:extends}
We have that $P_v$ is a $\Gamma_0(n)$--polygon;
furthermore, in the corresponding side pairing rule of $P_v$
 the $(2,1)$ component of each non-translation element is $n$.
\end{lem}
\bp
Let $(\bar a, \bar b)$ and $(\bar x, \bar y)$ be Farey pairs from $\Farey_v$ such that
\[
\bar a< \bar b,\quad \bar x<\bar y,\]
and such that $\bar a$ and $\bar b$ are non-adjacent in $\Farey_v$. Let us denote by $a,b,x,y$ the corresponding denominators, so that $a+b\le v$.
Lemma~\ref{lem:elem} implies that $n\not\divides ax+by$.
Using Lemma~\ref{lem:2-1-compo}, we see that the geodesics $(\bar a,\bar b)$ and $(\bar x,\bar y)$ are never identified, even up to orientations, under the action of $\Gamma_0(n)$; this implies that $P_v$ is a $\Gamma_0(n)$--polygon.

Suppose now that $N$ is the $(2,1)$ component of an element $g$ in the side-pairing rule that does not fix $\infty$.
Let $(a,b)$ be a non-free side $e$ of $dF$ 
so that $g_e=g$.
Note that $(a,b)\neq (0, 1)$ or $(1, 0)$.
If $e$ is paired or even, then we see from Lemma~\ref{lem:2-1-compo} that
some (possibly the same) side $(x,y)$ of $dF$ satisfies \[ n\divides N=ax+by.\]
Since $a,b,x,y \le \sqrt{n}$ and $\gcd(a,b)=1$,
we have that
\[
0<N=ax+by \le (a+b)\max(x,y)<2n,\]
and that $N=n$.
If $(a,b)$ is odd, we also see from Lemma~\ref{lem:2-1-compo} that
\[
n\divides N=a^2+ab+b^2<3n.\]
Using a simple parity argument, we conclude that $N=n$.
\ep

Let us record an algebraic consequence of our discussion so far.
\begin{lem}\label{lem:min-cusp}
For all integer $n\ge2$, we have that
\[
\Phi(n)\le \frac{[\Gamma\co\Gamma_0(n^2)]-v_3(n^2)}3.\]
When $n$ is a power of a prime $p$, we also have that
\[
 \Phi({\floor{\sqrt{n}}})  \le  \floor*{\frac{n(1+1/p)}3}.\]
\end{lem}
\bp
 Lemma~\ref{lem:extends} implies 
 \[
 \Phi(\floor{\sqrt{n}})=\#\cusp(P_v)-2\le u(n),\]
 from which the conclusion is immediate.
\ep

As a slight variation of the asymptotic formula used in the proof of Lemma~\ref{lem:lower1}, we have a (non-asymptotic) bound
\[
\Phi(m)\le 3m^2/\pi^2+2m\log m\]
 for all  $m\ge2$.
We can now establish the lower bound for $m(\Gamma_0(n))$.
\bp[Proof of part (1) in Theorem~\ref{t:bound}]
From 
Lemmas~\ref{lem:min-cusp} and \ref{lem:lower1}, 
we see that
$
\Phi(v)\le u(n)$, and hence, that
$m(\Gamma_0(n))\ge v$.

The equality $m(\Gamma_0(n))=v$ holds iff the polygon $P_v$ is maximal,
which is equivalent to that
$
u(n)
=\Phi(v)$.
As in the proof of Lemma~\ref{lem:lower1},
 this equality would imply
\[
\frac{n-\sqrt{n}}3<
u(n)
=\Phi\left(\floor{\sqrt{n}}\right)
<
 \frac{3}{\pi^2}n + \sqrt{n} \log n.\]
Note that the following function is strictly positive
for all $x>10^6$:
\[
f(x):=\frac{x-\sqrt{x}}3 - \frac{3}{\pi^2}x - \sqrt{x}\cdot {\log x}.\]
On the other hand, it is a routine computer verification to see
that the equation $u(n)
=\Phi(v)$ holds for an integer $n<10^6$ iff $n$ is in the given list. This completes the proof.
\ep

\subsection{Farey triple}
We now establish an efficient upper bound for $m(\Gamma_0(n))$ when $n=p,p^2$ or $pq$ for primes $p$ and $q$, with $p$ and $q$ being close.
A crucial idea is that we can group many (sometimes all) of the free sides of $dF$ into \emph{Farey $n$--triples}, 
each of which consists of three sides with certain properties that control the largest denominators in the resulting special polygons.
Recall our notation
\[
T:=\begin{pmatrix}1&1\\0&1\end{pmatrix}.\]
We emphasize our convention again that $n$ is only assumed to be a positive integer, unless further specified.

\bd\label{d:triple}
A sequence of distinct co-prime pairs of positive integers 
\[
\{(a_i,b_i)\co i=0,1,2\}
\]
is called an \emph{$n$--Farey triple}
if
the following holds for $i=0,1,2$ with the indices taken cyclically:
\begin{equation}\tag{F}\label{eqF}
\begin{pmatrix}
a_{i+1} & b_{i+1}
\end{pmatrix}
T
\begin{pmatrix}
a_{i} \\ b_{i}
\end{pmatrix}=n.
\end{equation}
\ed

\begin{figure}[htb!]
\leavevmode \SetLabels
\L(.1*.13) ${a_0}$\\%
\L(.15*.13) $a_0+b_0$\\%
\L(.17*.59) $\Delta_0$\\%
\L(.31*.13) $b_0$\\%
\L(.3*.69) $e_0$\\%
\L(.37*.13) $a_1$\\%
\L(.45*.13) $a_1+b_1$\\%
\L(.62*.13) $b_1$\\%
\L(.62*.69) $e_1$\\%
\L(.47*.69) $\Delta_1$\\%
\L(.66*.13) $a_2$\\%
\L(.75*.13) $a_2+b_2$\\%
\L(.88*.69) $e_2$\\%
\L(.88*.13) $b_2$\\%
\L(.8*.65) $\Delta_2$\\%
\endSetLabels
\begin{center}
\AffixLabels{\centerline{\includegraphics[width=10cm,angle=0]{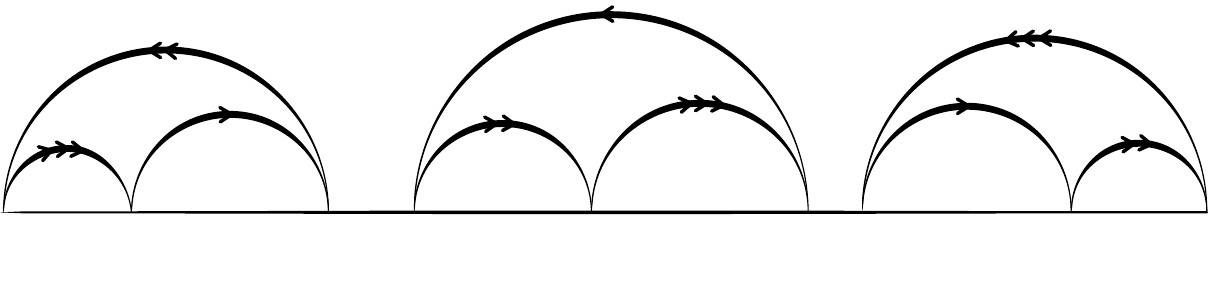}}}
\vspace{-24pt}
\end{center}
\caption{An $n$--Farey triple realized as denominators of Farey pairs. Single, double or triple arrows illustrate how these sides are identified by the action of $\Gamma_0(n)$. }
\label{fig:Pants1}
\end{figure}

Let $\{(a_i,b_i)\mid i=0,1,2\}$ be an $n$--Farey triple,
and pick the unique Farey pair $(\bar a_i,\bar b_i)$ corresponding to each $(a_i,b_i)$; see Lemma~\ref{lem:side}.
By Lemma~\ref{lem:2-1-compo}
and by the unique transitivity of the action of $\Gamma$ on special triangles, we have a unique $\gamma_i\in \Gamma_0(n)$ such that 

\begin{equation}\tag{\&}\label{eq:triple}
\left(
\bar a_i , 
\bar a_i \oplus \bar b_i,
 \bar b_i 
\right)
\stackrel{\gamma_i}{\to}
\left(
\bar a_{i+1}  \oplus \bar b_{i+1},
\bar b_{i+1}  ,
\bar a_{i+1}\right)
\end{equation}
The identifications of these triangles are illustrated in Figure~\ref{fig:Pants1}
by single, double or triple arrows.

The following lemma explains the connection between the (algebraically defined) $n$--Farey triples
and the (geometrically defined) sides of $d\Farey_v^*$. It further implies that
the first two pairs in an $n$--Farey triple determine the remaining pair.
\begin{lem}\label{lem:012}
Suppose we have a sequence of co-prime pairs of positive integers
\[Z=\{(a_i,b_i)\co i=0,1,2\}.\]
\be
\item\label{p:triple}
If $Z$ is an $n$--Farey triple,
then we have that
$
a_i+b_i=b_{i+1}+a_{i+2}$
for each $i$, 
and that
$
\min_i(a_i+b_i) < \sqrt{4n/3}$.
If we further assume that 
$\min_i (a_i+b_i)>\sqrt{n}$, then $Z$ consists of free sides in $d\Farey_v^*$ satisfying $ \max(a_i,b_i)<\sqrt{n}$ for all $i$.
\item\label{p:012-2}
Conversely, if we have 
 that $a_1(a_0+b_0)+b_1b_0=n$
and that 
\[
(
a_2,b_2)=(
a_0+b_0-b_1,
-a_0+a_1+b_1),
\]
then either $Z=\{(a_0,b_0)\}$
or $Z$ is an $n$--Farey triple.
\item\label{p:012-triple}
If $Z$ consists of sides in $d\Farey_v^*$,
and if \[a_{i+1}(a_i+b_i)+b_{i+1}b_i=n\]
for $i=0,1$,
then  either $Z=\{(a_0,b_0)\}$
or
$Z$ is an $n$--Farey triple.
\ee
\end{lem}
\bp
(\ref{p:triple})
We note that
\[
\begin{pmatrix}
n \\ n
\end{pmatrix}
= 
\begin{pmatrix}
\begin{pmatrix}
a_1 & b_1
\end{pmatrix}
T^t
\\
\begin{pmatrix}
a_0&b_0
\end{pmatrix}
T
\end{pmatrix}
\begin{pmatrix}
a_2\\b_2
\end{pmatrix}=
\begin{pmatrix}
a_1+b_1 & b_1 \\ a_0 & a_0+b_0
\end{pmatrix}
\begin{pmatrix}
a_2\\b_2
\end{pmatrix}.
\]
The matrix in the right-most term has the determinant $n$ and hence,
\[
\begin{pmatrix}
a_2\\b_2
\end{pmatrix}
=\begin{pmatrix}
a_0+b_0 & -b_1\\
-a_0 & a_1+b_1
\end{pmatrix}
\begin{pmatrix}
1 \\ 1
\end{pmatrix}
=
\begin{pmatrix}
a_0+b_0-b_1\\
-a_0+a_1+b_1
\end{pmatrix}.
\]
This proves the claimed equality. We also note
for each $i$ that
\[
a_ia_{i+1}+(a_i+b_i)b_{i+1}=n
    \Leftrightarrow (a_i+b_i)(a_{i+1}+b_{i+1})-n=b_ia_{i+1}.
\]
Setting $A=\min_i (a_i+b_i)$, we obtain that
\[
\left(1-\frac{n}{A^2}\right)^3\le
\prod_i\left(1-\frac{n}{(a_i+b_i)(a_{i+1}+b_{i+1})}\right)=\prod_i\frac{a_ib_i}{(a_i+b_i)^2}<\frac1{64}.\]
 This implies that $A< \sqrt{4n/3}$, as required.
 
For the latter conclusion we note for each $i$ that
\[
a_{i+1}\sqrt{n}<a_{i+1}(a_i+b_i)<a_{i+1}(a_i+b_i)+b_{i+1}b_i=n.\]
By a similar consideration for $b_i$, we have that $ \max(a_i,b_i)<\sqrt{n}$ for all $i$.
Since $\gcd(a_i,b_i)=1$, we see from Lemma~\ref{lem:side} that $d\Farey_v^*$ has $(a_i,b_i)$ as a side.

Assume for contradiction that $(a_0,b_0)$ is not a free side.
If $(a_0,b_0)$ is a paired side or an even side, then there exists a side $(a,b)$ in $dF$ such that $a_0a+b_0b=n$. 
In the case when $b>a$, we deduce from \[n=a(a_0+b_0)+(b-a)b_0=a_1(a_0+b_0)+b_1b_0\]
and from $0<a_1,b_1\le v<a_0+b_0$ that 
\[(a_1,b_1)=(a, b-a).\]
This contradicts $a_1+b_1>v$.
Similarly, if $b<a$, we would see from \[
n=b(a_0+b_0)+(a-b)a_0=b_2(a_0+b_0)+a_2a_0\]
that  $(a_2,b_2)=(a-b, b)$, which again is a contradiction.

If $(a_0,b_0)$ is an odd side, then we would have that \[n=a_0(a_0+b_0)+b_0^2=a_1(a_0+b_0)+b_1b_0,\]
and that $(a_1,b_1)=(a_0,b_0)$. This is again a contradiction, and we conclude that each $(a_i,b_i)$ is a free side.

(\ref{p:012-2})
Simply reading backward the former half of the proof of part (\ref{p:triple}),
we see that the equation~\eqref{eqF} in Definition~\ref{d:triple} holds for $i=0,1,2$. It is easy to see from the hypothesis that either $Z$ is a singleton or $Z$ consists of three distinct pairs.

(\ref{p:012-triple})
We have two solutions $(a_0+b_0,b_0)$ and $(a_2,a_2+b_2)$ of the equation
\[
a_1x+b_1y=n.\]
Using the hypothesis that $(a_i,b_i)$ is a side, we see 
\[ (a_2,a_2+b_2)=(a_0+b_0-b_1,b_0+a_1).\]
Then we compute
\[
(a_2+b_2)a_0+b_2b_0
=
(b_0+a_1)a_0+(a_1+b_1-a_0)b_0=n.\]
The conclusion now follows from part~(\ref{p:012-2}).
\ep

The following may be considered as the key strategy for the proof of Theorem~\ref{thm:bound}.
We say a set $A\sse\bZ^2$ is said to be \emph{primitive} if each vector $(x,y)\in A$ satisfies that $\gcd(x,y)=1$.
Let us introduce the notation
\[
S(a,b;n):=\{(x,y)\in\bZ^2\mid ax+by=n\}.\]

\begin{lem}\label{lem:primitive}
If $(a,b)$ is a free side of $d\Farey_v^*$,
and if the set
\[S(a,b;n)\cap\left((v-b,v]\times\bZ_+\cup \bZ_+\times(v-a,v]\right)\] 
is primitive, then $(a,b)$ belongs to an $n$--Farey triple, which consists of free sides in $d\Farey_v^*$.
\end{lem}
\bp
Let us establish a series of claims for the proof.
\begin{claim}\label{cla:triple1}
If $(x,y)\in S(a,b;n)\cap(0,v]^2$, then $x+y>v$ and $\gcd(x,y)\ne1$.
\end{claim}
Let $(x,y)$ be as in the hypothesis. We note from Lemma~\ref{lem:elem} that 
$x+y>v$. If $\gcd(x,y)=1$ then $(x,y)$ would be a side of $dF$, which contradicts the assumption that $(a,b)$ is free; see Lemma~\ref{lem:2-1-compo}. The claim follows.

\begin{claim}
We have that $\max(a,b)<\sqrt{n}$.
\end{claim}
The claim is vacuous if $v<\sqrt{n}$. So, let us assume $a=v=\sqrt{n}$, so that  $b<v$,
and that $(a-b,a)\in A\cap(0,v]^2$.
This contradicts Claim 1, since  $\gcd(a-b,a)=\gcd(a,b)=1$.
By symmetry, we see $b\ne \sqrt{n}$ as well.

Let us now consider the unique pair $(x_2,y_2)\in S(a,b;n)$ such that 
\[
0\le v-b< x_2\le v.\]
Since $x_2+b, a+b\in [v+1,\infty) \sse (\sqrt{n},\infty)$, we have
\[
0< n-av \le n-ax_2 =by_2<\sqrt{n}(x_2+b)-ax_2<b(x_2+\sqrt{n}).\]
In particular, we have 
$0<y_2<x_2+\sqrt{n}$.
We see from the hypothesis that
$\gcd(x_2,y_2)=1$.
 Claim~\ref{cla:triple1} implies $y_2>v\ge x_2$. 
 
Similarly, the unique $(x_0,y_0)\in S(a,b;n)$ satisfying \[0\le v-a<y_0\le v,\]
will have the properties that 
$\gcd(x_0,y_0)=1$
and that  \[y_0\le v<x_0<y_0+\sqrt{n}.\]
Then $(a_1,b_1):=(a,b)$, $ (a_2,b_2):=(x_2,y_2-x_2)$ and $(a_0,b_0):=(x_0-y_0,y_0)$ satisfy the equation~\eqref{eqF}
in Definition~\ref{d:triple}
with $i=0,1$. Furthermore, 
$Z:=\{(a_i,b_i)\}_{i=0,1,2}$ is not a singleton since $(a,b)$ is not odd. 
Lemma~\ref{lem:012} (\ref{p:012-triple}) implies that 
$Z$ is an $n$--Farey triple
consisting of free sides.\ep

Let us now consider the special case that $n$ is a prime or its square.
\begin{lem}\label{lem:farey-triple} 
If $n$ is a prime or the square of a prime, 
and if $(a,b)$ is a free side in $d \Farey_v^*$,
then the solution set $S(a,b;n)\cap\bZ_+^2$ is primitive. \end{lem}
\bp 
Let $(x,y)\in\bZ_+^2$ satisfy $ax+by=n$.
We
 have that
\[\gcd(x,y)\le \min(x,y)\le n/(a+b)<\sqrt{n}.\]
Since $n$ is a prime or its square,
we obtain that $\gcd(x,y)=1$.\ep

Recall $\|\cdot\|_F$ denotes the Frobenius norm.
We note the following simple fact on Möbius maps moving geodesics in the strip $0\le\im z\le 1$.
\begin{lem}\label{lem:strip-invert}
Let us consider a matrix
\[
g=\begin{pmatrix} a & b\\ c& d\end{pmatrix}\in\SL(2,\bZ)\]
with $c>0$.
If there exist some rational numbers $0\le x<y\le 1$ satisfying
\[
0\le g(y)<g(x)\le 1,
\]
then we have that
$|\tr g|\le c-2$ and that
$\|g\|_F<2c-1$.
\end{lem}
\bp
Using the equality $ad-bc=1$
and the sequence of inequalities
\[
0\le g(y)\le g(1)<g(\infty)<g(0)\le g(x)\le 1,
\]
we can establish
\[
-c<d\le b<0<-b\le a<c.\]
The conclusion is then immediate.
\ep

\bp[Completing the proof of Theorem~\ref{thm:bound}]
It is clear from the preceding lemmas that every free side of $dF$ belongs to a unique $n$--Farey triple.
We can hence partition the free sides of $dF$ into a collection of $n$--Farey triples
\[
\bigsqcup_{i=1,\ldots,k} \{(a^i_j,b^i_j)\co j=0,1,2\}\]
such that for all $i$ we have
\[a^i_0+b^i_0=\min_j (a^i_j+b^i_j)\in\left(\sqrt{n},\sqrt{4n/3}\right).\]
There exists a Farey pair $(\bar a^i,\bar b^i)$ corresponding to the free side $(a_0^i,b_0^i)$.
Set
\[
c_i:=\bar a^i \oplus \bar b^i \in \Farey^*_{\floor{\sqrt{4n/3}}}.\]
Then we have proved above that the newly added sides $(\bar a^i ,c_i)$ and $(c_i,\bar b^i)$ are paired to two of the free sides in $\Farey_v^*$.
In this way, we obtain the side-pairing rule for 
\[
P':=\hull(\Farey_v^*\sqcup \{c_i\}_i),\]
which is a maximal $\Gamma_0(n)$--polygon. In particular, we have
\[
m(\Gamma_0(n))\le
\max_{i=1,\ldots,k}\min_{j=0,1,2}(a^i_j+b^i_j)\le\floor*{\sqrt{4n/3}}.\]

Combining the definition of a Farey $n$--triple with Lemma~\ref{lem:extends}, 
we can choose a side pairing rule for $P'$ so that the $(2,1)$ component of each side-pairing element is $n$, with the sole exception of the translation $(x,y)\mapsto (x\pm1,y)$.
From Lemma~\ref{lem:strip-invert},  we see 
each non-translational side pairing element satisfies the claimed bounds on the trace
and the Frobenius norm.
\ep

Let us conclude this section by establishing an analogous result to Theorems~\ref{t:main} and \ref{t:bound}
 when $n$ is the multiplication of two sufficiently close odd primes.
 This theorem covers, for instance, all twin primes.

\begin{thm}\label{thm:twin}
Suppose we have two odd primes $p$ and $q$ satisfying \[\sqrt{p}<\sqrt{q}<\sqrt{p}+\sqrt{2}.\]
Then for $n:=pq$, the group $\Gamma_0(n)$ admits 
a free product decomposition
\[
\Gamma_0(n) =\bigast_{i=1}^\ell\form*{g_i}\bigast \form{T }\]
for some  matrices $g_i$ such that 
\[
(g_i)_{2,1}=
\begin{cases}
2n, &\text{ if }i\le q-p,\\
n, &\text{ if }i>q-p.\end{cases}
\]
We can further require that 
$|\tr g_i|\le (g_i)_{2,1}-2$ 
and $\|g_i\|_F<2(g_i)_{2,1}-1$.
Moreover, we have that
\[
\floor*{\sqrt{n}}\le m(\Gamma_0(n))\le\max\left(
\floor*{\sqrt{4n/3}},q\right).\]
\end{thm}

With the extra assumption that $q>83$, 
we have that $\floor*{\sqrt{4n/3}}\ge q$ above, 
and hence obtain the same bound as in Theorem~\ref{t:bound}:
\[
\floor*{\sqrt{n}}\le m(\Gamma_0(n))\le
\floor*{\sqrt{4n/3}}.\]
\bp[Proof]
Set $2k=q-p$.
The given bound on $\sqrt q - \sqrt p$ implies $v=p+k-1$.
We note that each element in the set
\[
A:=\{(k,p),(p,k)\}\cup\{(i,q-i)\mid k<i<p+k\}\]
is a free side of $P_v$.
Furthermore, if
$(a,b)$ is a free side  of $P_v$  not belonging to $A$,
then it is not hard to see that 
the set 
\[S(a,b;n)\cap((v-b,v]\times\bZ_+\cup \bZ_+\times(v-a,v])\] is primitive.
Applying Lemma~\ref{lem:primitive}
and suitably splitting the sides in $A$,
we obtain another $\Gamma_0(n)$--polygon $P'$
each of whose sides is exactly one of the following types:
\be[(i)]
\item a non-free side of $P_v$;
\item a pair that belongs to an $n$--Farey triple;
\item a side in the set $A':=\{(p+k,p),(p,k)\}$;
\item a side in the set
$
A'':=\{(i,q),(q,p-i)\mid k\le i\le p-k-1\}$;
\item a side in the set
$
A''':=\{(p-k+i,q),(q,p+k-i)\mid 0\le i\le 2k-1\}$.
\ee
The sides in $A'$, $A''$ and $A'''$ can be paired among themselves respectively; concretely, we have
\begin{align*}
(p+k,p)\cdot(p,k)&=n.\\
(i,q)\cdot(q,p-i)&=n.\\
(p-k+i,q)\cdot(q,p+k-i)&=2n.
\end{align*}
In particular, $P'$ is a maximal $\Gamma_0(n)$--polygon.
Using a similar argument to the case of $n=p$ or $p^2$, we obtain the conclusion.
\ep

\begin{rem}\label{rem:estimate}
\be
\item
For a positive integer $n$, 
let us denote by $k(n)$ the number of the $n$--Farey triples among the (necessarily free) sides in $d\Farey_v^*$. Algebraically, this number coincides with the number of ordered pairs $(a_0,b_0)$ of co-prime positive integers such that
for some choices of positive integers $a_1,b_1,a_2$ and $b_2$, all of the following hold:
\be[(i)]
\item 
For each $i$ with the indices taken cyclically, we have:
\[\begin{pmatrix}
a_{i+1} & b_{i+1}
\end{pmatrix}
T
\begin{pmatrix}
a_{i} &  b_{i}
\end{pmatrix}^t=n;
\]
\item $(a_i,b_i)$ is a co-prime pair for each $i$;
\item 
We have
\[
\sqrt{n}<a_0+b_0=\min_i(a_i+b_i)<a_1+b_1.
\]
\ee
The condition (iii) guarantees that one Farey $n$--triple is not counted more than once. Note we allow $a_0+b_0=a_2+b_2$.
\item
Let $n$ be a positive integer,
and let $\pi$ denote the quotient map from $ \bH^2$ onto the orbifold $\bH^2/\Gamma_0(n)$.
We define $P_v'$ as the union of $P_v=\hull(\Farey_v^*)$ with one special triangle glued on each odd side.
In the case when $n$ is a prime or its square, our proof of Theorem~\ref{thm:bound} implies the existence of one-to-one correspondence between the set of the connected components in $\left(\bH^2/\Gamma_0(n)\right)\setminus  \pi(P_v')$
and the set of Farey $n$--triples $\{(a_i,b_i)\}_i$ satisfying $a_i+b_i>\sqrt{n}$ modulo the cyclic $\bZ/3\bZ$--symmetry.
As a result, 
we have an estimate of the Euler totient summatory function
\[
\sum_{i=1}^{\floor{\sqrt{n}}} \phi(i)= u(n)-k(n)=\floor*{\frac{n(1+1/p)}3}-k(n).\]
\item In the case when $n=pq$ with $p$ and $q$ satisfying the hypotheses of Theorem~\ref{thm:twin}, each connected component of $\left(\bH^2/\Gamma_0(n)\right)\setminus \pi(P_v')$ is either a Farey triangle corresponding to some Farey $n$--triple,
or an ideal polygon with a puncture inside;
see Figure~\ref{fig:estimate}
for the description of this ideal polygon when $q=p+2$.
Such an ideal polygon has $p+1$ Farey triangles. Hence, we have an estimate
\[
\sum_{i=1}^{\floor{\sqrt{n}}} \phi(i) + k(n)+p+1=u(n)=\frac{
(p+1)(q+1)-v_3(pq)
}3.\]
\ee
\end{rem}

\begin{figure}[htb!]
  \tikzstyle {a}=[red,postaction=decorate,decoration={%
    markings,%
    mark=at position .5 with {\arrow[red]{stealth};}}]
  \tikzstyle {b}=[blue,postaction=decorate,decoration={%
    markings,%
    mark=at position .43 with {\arrow[blue]{stealth};},%
    mark=at position .57 with {\arrow[blue]{stealth};}}]
  \tikzstyle {v}=[draw,circle,fill=black,inner sep=0pt]
 \begin{tikzpicture}[thick,scale=1.2]
  \foreach \i in {2,...,6} {
    \draw  (360/10*\i:1.5)--(360/10*\i+360/10:1.5);
    \node [] at (360/10*\i+18:.5) {\tiny $q$};
    }
  \foreach \i in {8,...,10}{
    \draw  (360/10*\i:1.5)--(360/10*\i+360/10:1.5);
    \node [] at (360/10*\i+18:.5) {\tiny $q$};
}
\node [] at (360/10*4.8:1.3) {\tiny $1$};
\node [] at (360/10*4.35:1.08) {\tiny $p+1$};

\node [] at (360/10*5.8:1.3) {\tiny $3$};
\node [] at (360/10*5.35:1.08) {\tiny $p-1$};

\node [] at (360/10*6.8:1.3) {\tiny $5$};
\node [] at (360/10*6.35:1.08) {\tiny $p-3$};

\node [] at (360/10*8.8:1.3) {\tiny $p$};
\node [] at (360/10*8.2:1.3) {\tiny $2$};

\node [] at (360/10*9.8:1.3) {\tiny $2$};
\node [] at (360/10*9.2:1.3) {\tiny $p$};

\node [] at (360/10*10.8:1.3) {\tiny $4$};
\node [] at (360/10*10.35:1.08) {\tiny $p-2$};

\node [] at (360/10*12.6:1.3) {\tiny $p-3$};
\node [] at (360/10*12.2:1.1) {\tiny $5$};

\node [] at (360/10*13.6:1.2) {\tiny $p-1$};
\node [] at (360/10*13.2:1.3) {\tiny $3$};
    \draw  (360/10*5+90:1) node [] {\tiny $\cdots$};
    \draw  (360/10*9+92:1) node [] {\tiny $\ddots$};
\foreach \i in {0,...,9} {
  \draw (0,0) -- (360/10*\i:1.5) node [v] {};
}
\draw (360/10*4:1.5) -- (360/10*4.5:3) node [v]  {};
\draw (360/10*5:1.5) -- (360/10*4.5:3);
\node [] at (360/10*4.5:2.5) {\tiny $p$};
\node [above] at (360/10*5:1.55) {\tiny $1$};
\node [below] at (360/10*4.15:1.85) {\tiny $p+1$};
\end{tikzpicture}
\caption{A complementary region of $\hull \left(\Farey_v^*\right)$ in $\bH^2/\Gamma_0(pq)$ with $q=p+2$. Three corners of each Farey triangle are labeled with the denominators of  corresponding cusp points.}
\label{fig:estimate}
\end{figure}

\section{Optimality of the upper bound}\label{s:optimal}
Let us say an integer $n$ is \emph{cashew} if for some $n$--Farey triple
$\{(a_i,b_i)\}$,
\[\min_i (a_i+b_i)=\floor*{\sqrt{4n/3}}.\]
The goal of this section is to establish the following, which implies the optimality of the coefficient $\sqrt{4/3}$ in Theorem~\ref{thm:bound}.
Note that for $n\ge37$ we have $\floor*{\sqrt{4n/3}}>{\sqrt{n}}$.

\begin{thm}\label{thm:upperbound1}
Let $n\ge37$ be a prime or its square.
Then we have $m(\Gamma_0(n))=\floor*{\sqrt{4n/3}}$ if and only if $n$ is a cashew integer.
\end{thm}

\begin{rem}\label{rem:cashew}
\be
\item By the proof that we give below, the same conclusion holds for $n=pq$ with $p$ and $q$ as in Theorem~\ref{thm:twin}, and with the extra assumption that $q>83$.
\item
One may drop the hypothesis $n\ge37$ in the above, although the proof will then require (rather light) computer verification. In particular, the following is the list of all cashew primes up to $500$:\\
5, 11, 17, 23, 41, 43, 47, 59, 71, 73, 89, 97, 101, 103, 107, 137, 139, 191, 211, 229, 233, 239, 241, 269, 281, 353, 389, 409, 419, 421, 431, 457, 499.\\
Notice that $5,11$ and $17$ also occur in the list of Theorem \ref{t:bound}.
\item One can verify by computer that for any cashew prime $p<150000$, the number of  $p$-Farey triples $(a_0,b_0,a_1,b_1,a_2,b_2)$ verifying the cashew property is  bounded by $2$ up to permutation and reflection. This property possibly holds for all cashew primes.
\ee
\end{rem}

We will utilize an alternative formulation of cashew integers as below. 

\begin{lem}\label{l:cashew2}
An integer $n\ge37$ is cashew iff 
there exist $s,t,a,b\in\bN_+$ such that
\begin{equation}\tag{C}\label{eqC}
\begin{cases}
n = sa+tb& \\
    s+t>a=\floor{\sqrt{4n/3}}>
t\geq 
b \geq a-s\ge1  & 
\end{cases}
\end{equation}
\end{lem}

\bp
Assuming $n$ is cashew, let us pick
a corresponding $n$--Farey triple
 $\{(a_i,b_i)\}$.
By Lemma~\ref{lem:012},
we may assume that
\[\sqrt{n}<\floor*{\sqrt{4n/3}}=
a_0+b_0=\min_i(a_i+b_i)<\max_i(a_i+b_i)=a_1+b_1,
\]
possibly after 
applying a cyclic permutation to $\{(a_i,b_i)\}$,
and possibly after changing the roles of $a_i$'s and $b_i$'s.
Lemma~\ref{lem:012}
implies that
the collection $\{(a_i,b_i)\co i=0,1,2\}$ consists of free sides in $d\Farey_v^*$ with $v:={\floor{\sqrt{n}}}$ 
such that $\max_i\max(a_i,b_i)<\sqrt{n}$.

Setting $(s,t,a,b):=(a_1, b_1, a_0+b_0,b_0)$, we have the following.
\begin{itemize}
\item $s+t-a = a_1+b_1-a_0-b_0>0$;
\item $a-t=a_0+b_0-b_1=a_2>0$;
\item $t-b = b_1-b_0=(a_1+b_1)-(b_0+a_1)=(a_1+b_1)-(a_2+b_2)\ge0$;
\item $b-(a-s)=b_0-(a_0+b_0)+a_1=(a_2+b_2)-(a_0+b_0)\ge0$;
\item $a-s=a_0+b_0-a_1> \floor*{\sqrt{4n/3}}-\sqrt{n}>0$;
    \item $a=\floor*{\sqrt{4n/3}}\leq \sqrt{4n/3}= \sqrt{4(sa+tb)/3}<a+1$.
\end{itemize}
This verifies the condition~\eqref{eqC}.

    Conversely, suppose $(s,t,a,b)$ satisfies \eqref{eqC}. 
We can uniquely choose a triple $\{(a_i,b_i)\}$ of positive pairs such that
\[
(s,t,a,b)=(a_1, b_1, a_0+b_0,b_0)\]
and such that
\[
(
a_2,b_2)=(a-t,s+t-(a-b))
=(
a_0+b_0-b_1,
-a_0+a_1+b_1).
\]
By Lemma~\ref{lem:012} (\ref{p:012-2}), 
the triple is an $n$--Farey triple.
It is easy to see that $\min_i(a_i+b_i)=a$, certifying that  $n$ is cashew.\ep

We can now establish that cashew primes or cashew prime squares $n$ realize the upper bound $\floor{\sqrt{4n/3}}$ of $m(\Gamma_0(n))$.
\bp[Proof of Theorem~\ref{thm:upperbound1}]
We have noted that $\floor*{\sqrt{4n/3}}>\sqrt{n}\ge v:=\floor*{\sqrt{n}}$.

($\Rightarrow$)  By hypothesis, we have $m(\Gamma_0(n))>v$.
Lemmas~\ref{lem:primitive} and~\ref{lem:farey-triple} imply that 
the (non-empty) collection of the free sides in $d\Farey_v^*$ can be partitioned into a disjoint collection of $n$--Farey triples \[
\sqcup_{1\le i\le k}\{(a^i_j,b^i_j)\co j=0,1,2\}.\]
By Lemma~\ref{lem:012}, we have that
\[
\floor*{\sqrt{4n/3}}=m(\Gamma_0(n))\le \max_i \min_j (a^i_j+b^i_j)
\le \floor*{\sqrt{4n/3}}.\]
This implies that $n$ is cashew, certified by some of the above $n$--Farey triples.

($\Leftarrow$) Let $s,t,a,b \in \bN_+$ be as in the equation \eqref{eqC} for $n>17$.
We have a maximal $\Gamma_0(n)$--polygon $P$ such that 
the denominator sequence $dX$ of the cusp set of $P$ satisfies $\max dX=m(\Gamma_0(n))$.
In the proof of Lemma~\ref{l:cashew2},
we have constructed a Farey $n$--triple $\{(a_i,b_i)\}$ 
containing $(s,t)$ such that \[s+t>a=\min_i (a_i+b_i)=\floor*{\sqrt{4n/3}}>\sqrt{n}.\]
By Lemma~\ref{lem:012}, this triple corresponds to free sides in $\Farey_v^*$ with $v:=\floor{\sqrt{n}}$ such that $\max_i\max(a_i,b_i)<\sqrt{n}$.
In particular, neither $s^2+t^2$ nor $s^2+st+t^2$ are multiples of $n$.

\begin{claim}\label{cla:mn-side}
We have that $(s,t)$ is a side of $dX$.
\end{claim}
Let $(\bar s,\bar t)$ be the Farey pair
corresponding to $(s,t)$ by Lemma~\ref{lem:side}.
Assume for contradiction that the geodesic $\gamma:=(\bar s,\bar t)\sse\bH^2$ is not a side of $P$. 
Since $s+t> m(\Gamma_0(n))$,
 the exterior of $P$
 contains the point $\bar s\oplus\bar t$.
It follows that $\gamma$ is contained in the exterior of $P$, 
possibly except for the endpoints. 
We can find a side $(\bar c,\bar d)$ of $P$ enclosing this geodesic $\gamma$ in $\bH^2$ along with the $x$--axis such that the corresponding side $(c,d)$ of $dX$ satisfies 
\[c \leq s,\qquad d \leq t,\qquad c+d \leq \max(s,t)\le v.\]
Since $(c+d)^2\le n$, we see that $(c,d)$ is neither even nor odd,
and that another side $(\alpha,\beta)$ of $dX$ is paired to $(c,d)$ such that $\alpha,\beta\le\max dX\le a$.
For some $k\ge1$, we see from the conditions in Lemma~\ref{l:cashew2} that
\[k(sa+tb)=kn=\alpha c+\beta d
\le a\cdot\max(s,t)
\le \max(sa,t(b+s))<sa+tb=n.\] 
This is a contradiction, and so Claim~\ref{cla:mn-side} follows.

The last part of the proof is the following.
\begin{claim}\label{cla:ab-side}
We have that either $(a,b)$ or $(a-t,a)$ is a side of $dX$.
\end{claim}
Suppose a side $(\alpha,\beta)$ of $dX$ is paired to $(s,t)$.
For some $k\geq 1$ we have
\[n\leq kn = \alpha s+\beta t \leq {2\max(s,t)a}< 2n.\]
Therefore, $k=1$ and $s\alpha+t\beta=sa+tb$. This equation has solutions of the form $(\alpha,\beta)=\left(a+tq,b-sq\right)$, where $q\in \mathbb{Z}$. 
In the case when $s>a-b$, the equation has a unique solution $(\alpha,\beta)=(a,b)$ satisfying $0<\alpha,\beta\leq a$. 
This implies that $(a,b)$ is a side in $dX$ paired with $(s,t)$.
If  $s=a-b$, then the equation has two solutions $(\alpha,\beta)=(a,b)$ and $(\alpha,\beta)=(a-t,b+s=a)$ satisfying $0<\alpha,\beta\leq a$;
in this case $(a-t,a)$ or $(a,b)$ is a side paired with $(s,t)$. Claim~\ref{cla:ab-side} is proved.

The above claim immediately implies the conclusion since
\[ \floor*{\sqrt{4n/3}}=a\le \max dX = m(\Gamma_0(n))\le \floor*{\sqrt{4n/3}}.\qedhere\]
\ep

\section*{Acknowledgements}
The authors are very grateful to Asbjørn Nordentoft for his interest in this work and for pointing us to the paper~\cite{Steiner2022}. The authors are also thankful to Yichen Tao for his help with some coding of initial versions of the algorithm for finding the special polygons. The first named author is supported by the Institute of Mathematics, Vietnam Academy of Science and Technology under the code IM-VAST01-2022.01.
The second named author is supported by Samsung Science and Technology Foundation under Project Number SSTF-BA1301-51. 
 The last named author is supported by the National University of Singapore academic research grant  A-8000989-00-00.

\bibliographystyle{amsplain}
\bibliography{ref}
\end{document}